\documentclass{amsart}
\theoremstyle{plain}
\newtheorem{Thm}{Theorem}
\newtheorem{Cor}{Corollary} 
\newtheorem{Lem}[Thm]{Lemma}

\numberwithin{equation}{section}

\newcommand{\CC}{{\mathbb{C}}}
\newcommand{\RR}{{\mathbb{R}}}

\newcommand{\NN}{{\mathbb{N}}}
\begin{document}

\title[Solving the Gleason problem on linearly convex domains]
{Solving the Gleason problem on linearly convex domains}
\author{Oscar Lemmers}
\email[Oscar Lemmers]{lemmers@science.uva.nl}
\author{Jan Wiegerinck}
\email[Jan Wiegerinck]{janwieg@science.uva.nl}
\date{June 27, 2001}

\begin {abstract}
Let $\Omega$ be a bounded, connected linearly convex set in $\CC ^n$ with 
$C^{1+\epsilon}$-
boundary. We show that the maximal ideal (both in $A(\Omega)$ and 
$H^{\infty}(\Omega)$) consisting of all functions vanishing at $p \in \Omega$ 
is generated by the coordinate functions $z_1 - p_1, \ldots, z_n - p_n$.
\end{abstract}
\maketitle

\section{Introduction}

\noindent Let $\Omega$ be a bounded domain in $\CC^n$. 
Let $R(\Omega)$ (usually $A(\Omega)$ or $H^{\infty}(\Omega)$) be a ring of 
holomorphic functions that contains the 
polynomials, and let $p = (p_1, \ldots, p_n)$ be a point in $\Omega$.
Recall the Gleason problem, cf. \cite{glea}:
is the maximal ideal in $R(\Omega)$ consisting of functions vanishing at $p$,
generated by the coordinate functions $(z_1 - p_1)$, $\ldots$, 
$(z_n - p_n)$ ? \\ 
One says that a domain $\Omega$ has the {\em Gleason $R$-property} 
if this is the case for all points $p \in \Omega$. We also say that it has 
the Gleason-property with respect to $R(\Omega)$.\\
Leibenzon was the first to solve a non trivial Gleason problem. He
proved (\cite{khen}) that the Gleason problem can be solved on any convex 
domain in $\CC^n$ having a $C^2$-boundary. This result was sharpened by 
Grang\'{e} (\cite{gran}, for $H^{\infty}(\Omega)$), and by Backlund and 
F\"{a}llstr\"{o}m (\cite {bafa1} and \cite{bafa2}, for $H^{\infty}(\Omega)$ 
and $A(\Omega)$ respectively), 
for convex domains in $\CC^n$ having only a $C^{1 + \epsilon}$-boundary.\\
Using his theorem on solvability of the $\overline{\partial}$-problem 
(\cite{oev1}), 
\O vrelid proved in \cite{oev2} that a strictly pseudoconvex domain 
in $\CC^n$ with $C^{2}$-boundary has the Gleason $A$-property. Forn\ae ss and 
\O vrelid showed in \cite{foev} that a pseudoconvex domain in $\CC^2$
with real analytic boundary has the Gleason $A$-property. This was extended
by Noell (\cite{noel}) to pseudoconvex domains in $\CC^2$ having 
a boundary of finite type.\\
Backlund and F\"{a}llstr\"{o}m proved in 
\cite{bafa4} that a bounded, pseudoconvex 
Reinhardt domain in $\CC^2$ with $C^{2}$-boundary that contains the 
origin, has the Gleason $A$-property. The present authors showed in \cite{LeWi}
that one does not need that $\Omega$ is pseudoconvex or that it contains the 
origin. They also solved the $H^{\infty}$-problem for such Reinhardt domains.
\\Note that there are not always solutions to the Gleason problem; in fact, 
Backlund and 
F\"{a}llstr\"{o}m showed (\cite{bafa3}) that there even exists an 
$H^{\infty}$-domain of holomorphy on which the problem is not solvable.\\
In this article, we return to the original method of Leibenzon, and use
it to solve the Gleason problem on $\CC$-convex 
domains (these are domains such that their intersection with any complex line 
passing through the domain is connected and simply connected) in $\CC ^n$ with 
$C^{1 + \epsilon}$-boundary. We denote the derivate of a function $g$ with 
respect to the $k$'th coordinate with $D_k g$. After translation we can 
assume that the domain contains the origin, and that $p=0$.
If $\Omega$ is convex, it is easy to see that 
for $f_i(z):= \int _{0} ^ {1} D_i f(\lambda z) d \lambda$, $f(z) = 
\sum_{i=1}^{n}z_if_i(z)$. The hard part is 
to show that $f_i \in A(\Omega)$. Leibenzon did this
by making estimates of $D_i f(\lambda z)$ on the
line segment between $0$ and $1$. If one considers $\CC$-convex domains, this
method doesn't work, of course. However, for a polynomial $P$, and
\[T_i(P)(z) := \int _{0} ^ {1} D_i P(\lambda z) d \lambda
,\] we still have that $P(z) = \sum_{i=1}^{n}z_iT_i(P)(z)$. The theorem of
Cauchy gives \[T_i(P)(z) = \int _{0} ^ {1} D_i P(\lambda z)
 d\lambda = \int _{\gamma_{z}} D_i P(\lambda z)
 d\lambda \] for any curve $\gamma_z$ in $\CC$ that connects
$0$ and $1$. We choose this curve $\gamma_z$ such that for all $s \in
[0,1]$ the point $\gamma_z(s)z$
is in $\overline{\Omega}$ intersected with the complex line through $0$ and
$z$. Estimating
$D_i P(\lambda z)$ on $\gamma_z$ yields an estimate
of $T_i(P)(z)$ in terms of $||P||_S$ (where $S$ is a suitable compact subset
of $\Omega$) instead of $||P||_{\tilde{\Omega}}$
(where $\tilde{\Omega}$ is the convex hull of $\Omega$). Then the fact that 
$\Omega$ is a Runge domain is used, first to extend $T_i$ to 
$H^{\infty}(\Omega)$, then to $A(\Omega)$.

\section{$\CC$-convex sets}
\noindent 
In $\RR ^{n}$ there are two natural definitions of convexity. A set $E$ is
convex if 
\begin{enumerate}
\item The intersection of $E$ with each line is connected, \quad \quad or
\item Through every point in the complement of $E$ there passes a hyperplane
which does not intersect $E$.
\end{enumerate}
If one assumes that $E$ is connected, these definitions are equivalent. These
definitions lead to the following generalizations of convexity in $\CC ^n$ :
\begin{enumerate}
\item One says that a set $\Omega \subset \CC ^n$ is $\CC${\bf -convex} if all 
its intersections with complex lines are connected and simply connected.
\item $\Omega$ is said to be {\bf linearly convex} (also : lineally convex) if
through every point in the complement of $\Omega$ there passes a complex 
hyperplane that does not intersect $\Omega$.
\item An open set $\Omega$ in $\CC ^n$ is called {\bf weakly linearly convex}
if through every point of $\partial \Omega$ there passes a complex 
hyperplane that does not intersect $\Omega$.
\end{enumerate}
For proofs of the following assertions and more information on 
$\CC$-convex sets we refer to \cite{APS1}, \cite{APS2} and \cite{Hor}.
\begin{itemize}
\item Every $\CC$-convex set is pseudoconvex. 
\item Every $\CC$-convex set is a Runge domain.
\item For a bounded connected domain $\Omega$ with $C^{1+\epsilon}$-boundary
all the previous definitions are equivalent, and every complex line passing 
through $\Omega$ will intersect $\partial \Omega$ transversally.
\end{itemize}

\section{Definitions and auxiliary results}
\noindent We
establish some notations : by $||f||_{\Omega}$ we denote the supremum of 
$|f|$ on $\Omega$. We denote the derivate of a function $g$ with 
respect to the $k$'th coordinate with $D_k g$.
The unique complex line through the points $0$ and 
$z$ is denoted by $L_{0,z}$. For a $w \in \partial \Omega$, we denote by 
$n_w$ the inner unit normal
vector to $\partial \Omega$ at $w$. Given $z \in \CC ^n$, we denote by
$\pi _{z} (n_w) = (\pi_{z}(n_w)_1, \ldots, \pi_{z}(n_w)_n)$ the orthogonal
projection of $n_{w}$ onto $L_{0,z}$. Constants may be denoted 
by the same letter at various places, even though their values are different.

\begin{Lem}
(Cf. \cite{bafa2}.)
Let $w \in \partial \Omega$. There exist a neighborhood $W$ of $w$ and a 
$\sigma
> 0$ such that for all $z \in W \cap \partial \Omega$, and $e$ a complex unit
tangent vector to $\partial \Omega$ at $z$, the following holds : if $0 <
s < 1$ then $|t| < (1- s) ^{1/(1 + \epsilon/2)} <  \sigma$ implies 
that $z + (1 - s)\pi _{z} (n_w) + te \in \Omega$.
\end{Lem}

\begin{proof}
Since $\Omega$ has $C^{1 + \epsilon}$-boundary, there is a $C^{1 + 
\epsilon}$-defining function $r$ with $\nabla r \neq 0$ on $\partial \Omega$, 
such that $\Omega = \{z \in \CC^n: r(z) < 0\}$. We compute at $z \in W \cap
\Omega$ :
\begin{eqnarray*}&& r(z + (1 - s)\pi _{z} (n_w) + te) - r(z) \\
&& = \Re 
\left(\sum_{j = 1} ^{n} \frac {\partial r} {\partial z_j} (z) ( 
(1 - s) \pi _{z} (n_w) _{j} + te_j) \right)+ 
O(|(1 - s)\pi _{z} (n_w)  + te|^{1 + \epsilon}) \\
&& =
 (1 - s) \Re \left(\sum_{j = 1} ^{n} \frac 
{\partial r} {\partial z_j} (z) \pi _{z} (n_w) _{j}\right) + 
(1 - s) g(s,t),\end{eqnarray*}
with
\[g(s,t) = O(|\frac{s -1}
{(1 - s)^{1/(1 + \epsilon)}} \pi _{z} (n_w) + \frac{te}{(1-s)
^{1/(1 + \epsilon)}}|^{1 + \epsilon}), \] since $e$ is a complex tangent 
vector to $\partial \Omega$ at
$z$, causing $\Re (\sum_{j = 1} ^{n} \frac {\partial r}{\partial z_j} (z) te_j)
$ to vanish. Let $1- s$, $t$ tend to zero with $|t| < (1- s)
^{1/(1 + \epsilon/2)}$. Then the term 
\[\frac{s -1}{(1 - s)^{1/(1 + \epsilon)}}\pi _{z} (n_w)  + 
\frac{te}{(1 - s)^{1/(1 + \epsilon)}}\] tends to zero. For every $1 - 
s$ close to zero, there is a small $t$ such that $z + (1- s)
\pi _{z} (n_w) + te \in \Omega$, in other words : $r(z + (1- s)\pi _{z} 
(n_w) + te) < 0$. Since $\pi _{z} (n_w)$ is not a complex tangent vector 
(if we choose $W$ small enough), $\Re (\sum_{j = 1} ^{n} \frac {\partial r}
{\partial z_j} (z) (1 - s) \pi _{z} (n_w) _{j}) \not = 0$. Hence this term 
has to be negative.
Then there is a constant $\sigma > 0$ such that $|t| < (1- s)
^{1/(1 + \epsilon/2)} < \sigma$ implies that \[(1 - s) \Re \left(\sum_
{j = 1} ^{n} \frac {\partial r} {\partial z_j} (z) \pi _{z} (n_w)_{j} 
\right) + (1 - s) g(s,t) < 0,\] in other words : $z + (1 - 
s)\pi _{z} (n_w) + te \in \Omega$. 
\end{proof}

\begin{Lem}
For every $z \in \partial \Omega$ there exists a neighborhood $W$ of $z$ and 
a $\sigma > 0$ such that $F_z(w,s) := w + 
(1-s)\pi_w(n_z)$ is a bijection from $W \cap \partial \Omega \times 
(1 - \sigma,1]$ to $W \cap \overline{\Omega}$.
\end{Lem}

\begin{proof}
Observe that $\pi_z(n_z) \not= 0$ (since $L_{0,z}$ intersects $\partial \Omega$
transversally), such that the Jacobian of $F_z(w,s)$ does not vanish at 
$(z,1)$. Now apply the inverse function theorem.
\end{proof}

\noindent Compactness of $\partial \Omega$ and lemma 1 and 2 now yield :

\begin{Cor}
There exist a finite number of open sets $W_1$, $\ldots$, $W_m$ in $\CC^n$ 
containing points $w_1$, $\ldots$, $w_m$ respectively, and a 
$\sigma > 0$ such that :
\begin{itemize}
\item $\partial \Omega \subset \cup W_i$
\item $F_{w_i}(w,s) : W_i \cap \partial \Omega \times (1 - \sigma,1] 
\rightarrow W_i \cap \overline{\Omega}$ is a bijection for all 
$1 \leq i \leq m$.
\item For $w \in W_i$, and $e$ a complex unit tangent vector to $\partial 
\Omega$ at $w$ : if $0 < s < 1$ and $|t| < (1- s) ^{1/(1 + 
\epsilon/2)} < \sigma$, then $w + (1 - s)\pi _{w} (n_{w_i}) + te \in \Omega
\quad \forall 1 \leq i \leq m$.
\end{itemize}
\end{Cor}

\noindent Let $V_i := \{z \in \overline{\Omega} : \exists w \in 
W_i, s \in (1 - \frac{1}{2}\sigma, 1] \; : z = w + (1 - s)\pi_{w}(n_{w_i})
\}$, $V := \cup _i V_i$. Since $\Omega \setminus V$ is compact in $\Omega$, 
we have that $A:=\min d(\partial \Omega, \Omega \setminus V) > 0$.

\vskip5mm \noindent {\bf Definition.}
For $z \in V$ we say that $w$, $w_k \in W_k$ {\bf correspond} to $z$ if there 
is an $s \in (1- \frac{1}{2} \sigma,1]$ such that $z = w + 
(1 - s)\pi_{w}(n_{w_k})$.\\
For $z \notin V$ we fix a smooth curve $\gamma_z$ in $\CC$, from $0$ to 
$1$, without loops, such that \[\gamma_z(s)z \in \Omega \; \;
\text{and} \quad d(\gamma_z(s)z, \partial \Omega) \geq A \quad \forall s \in 
[0,1].\]
For $z \in V$ with corresponding $w$ and $w_k$, we fix a smooth curve 
$\gamma_z$ in $\CC$, from $0$ to $1$, without loops,
that consists of two parts $\gamma_z ^1$ and $\gamma_z ^2$. 
We choose $\gamma_z^{2}$ such that for $s \in [1 - \frac{1}{2} \sigma,1]$,
$\gamma_z^{2}(s)z = z + (1-s)\pi_{w}
(n_{w_k})$, thus, with $\mu_z$ the unique constant such that 
$\mu_z z = \pi_{w}(n_{w_k})$ : $\gamma_z^{2}(s) = 
1 + (1-s)\mu_z$ for all $s \in [1 - \frac{1}{2}\sigma,1]$. Note that for $z \in
V_k$, $\mu_z$ is bounded.\\ We choose the curve $\gamma_z ^{1}$ in $\CC$, 
from $0$ to $1 - \frac{1}{2}\sigma$,
without loops, such that \[\gamma_z ^{1}(s)z \in \Omega \; \; \text{and} 
\quad d(\gamma_z^{1}(s)z, \partial \Omega) \geq A \quad \forall s \in 
[0, 1 - \frac{1}{2} \sigma].\] 
It is possible to choose the curves $\gamma_z$ such that there is a constant 
$M$ with $|\gamma^{'}_{z}(s)| \leq M \; \; \forall z$, $s$, and that
$z_n \rightarrow z$ implies that $\gamma_{z_n}'(s) \rightarrow 
\gamma_{z}'(s)\; \; \forall z$, $s$.

\vskip5mm
\noindent {\bf Definition.}
For a polynomial $P$ that vanishes at $0$ we define linear operators 
$T_i$ $(1 \leq i \leq n)$ as follows : 
\[ T_i(P)(z) := \int^{1}_{0} D_i P(\lambda z) d\lambda .\] 
The operators are clearly related to the differentiated simplex functionals 
$\mu_{a,b}(\partial f)= \int _{0}^{1} \partial f(a+t(b-a)) dt$, that are
studied in \cite{APS2}, 3.2.
Note that one has \[P(z)= \int^{1}_{0} \frac{dP(\lambda z)}{d \lambda}
 d\lambda = \int^{1}_{0} \sum_{i=1}^{n} z_i D_i P(\lambda z)
 d\lambda = \sum _{i=1}^{n} z_i T_i(P)(z) .\]

\section{The key estimate}

\begin{Lem} Let $P$ be a polynomial that vanishes at $0$.
Let $p \in \Omega$. 
There exist a constant $K$ that depends only on $\Omega$, a neighborhood $B$ of
$p$ and a compact set $S$ in $\Omega$, such that $||T_i(P)||_B \leq K||P||_S$.
\end{Lem}

\begin{proof}
First we consider the case that $p \notin V$. Then we choose $B \subset \Omega
$ such that 
$\overline{B} \cap V = \emptyset$. 
Let $z \in \overline{B}$. Then $d(\gamma_z(s)z, 
\partial \Omega) \geq A$ for all $s \in [0,1]$. If $a_i$ is the $i$'th 
unit vector, we have
\begin{eqnarray*}
T_i(P)(z)& =& \int^{1}_{0} D_i P(\lambda z)
d\lambda = \int_{\gamma _{z}} D_i P(\lambda z)
d\lambda\\ & =&\int_{\gamma _{z}} \frac{dP(\lambda z + ta_i)}{dt}|_{t=0} \; 
d\lambda = \int_{0}^{1} \frac{dP(\gamma(s) z + ta_i)}{dt}|_{t=0} 
\gamma^{'}(s) ds.
\end{eqnarray*}
We construct a compact $S$ such that for all $s \in [0,1]$ $S$ contains 
the circle in the complex line through $sz$ and $sz + a_i$ with center $sz$ 
and radius $A$. Then we have for all $s \in
[0,1]$ : \[\left|\frac{dP(\gamma(s)z + ta_i)}{dt}|_{t=0} \right| = 
\left|\frac{1}{2 \pi i} \int _{C(0,A)} \frac{P(\gamma(s)z + ta_i)}
{t^2}dt\right| \leq \frac{||P||_S}{A},\]
where $C(0,A)$ is the circle with center $0$ and radius $A$. Hence
\[\left|\int_{\gamma _{z}} \frac{dP(\lambda z + ta_i)}{dt}|_{t=0} d\lambda
\right| \leq \frac{M ||P||_S}{A} \leq K ||P||_S.\]
Now we consider the case that $p \in V \cap \Omega$. Choose $B$ such that 
$\overline{B} \subset V \cap \Omega$.
Let $z \in V \cap \overline{B}$. Take $w$, $w_k \in \partial \Omega$ 
corresponding to $z$. One can make the appropriate estimate on 
$\gamma_{z}^{1}$ as above. Let $e$ be a complex unit tangent vector to 
$\partial \Omega$ at $w$. We construct $S$ such that for all $s \in 
[1 - \frac{1}{2} \sigma ,1]$ $S$
contains the circles in the complex line through $z +
(1 - s)\pi_{w}(n_{w_k})$ and $z + (1 - s)\pi_{w}(n_{w_k}) + e$ with center 
$z + (1 - s)\pi_{w}(n_{w_k})$ and radius $s^{1/(1+\epsilon/2)}$. Then
\begin{eqnarray*}
&&\left|\frac{dP(\gamma(s) z + te)}{dt}|_{t=0} \right| = \left|\frac{dP(z + 
(1- s) \pi_{w}(n_{w_k}) + te)}{dt}|_{t=0} \right| \\ 
&&= \left|\frac{1}{2 \pi i} \int _{C(0,s^{1/(1+\epsilon/2)})}
\frac{P(z + (1 - s) \pi_{w}(n_{w_k}) + te)}{t^2}dt\right| \leq \frac{||P||_S}
{s^{1/(1+\epsilon/2)}},
\end{eqnarray*}
where $C(0,s^{1/(1+\epsilon/2)})$ is 
the circle with center $0$ and radius $s^{1/(1+\epsilon/2)}$. Hence
\[\left|\int _{\gamma_z ^2} \frac{dP(\lambda z + te)}{dt}|_{t=0} \; d \lambda 
\right| \leq \int ^{1}_{0}
\frac{M||P||_S}{s^{1/(1+\epsilon/2)}} d s \leq K ||P||_S.\]
We can choose $S$ to be compact in $\Omega$.

\vskip5mm \noindent
For $z \in V_k$ the corresponding $w$ depends continuously on $z$, hence we
can choose linearly independent complex unit tangent vectors $e^1(z)$, 
$\ldots$, $e^{n-1}(z)$ to $\partial \Omega$ at $w$ that depend continuously 
on $z$.
As a consequence of the theorem of Cauchy and the chain rule, we have :
\[\int_{\gamma_z} \frac{dP(\lambda z + te^j(z))}{dt} |_{t=0} \; d\lambda =
\int_{0} ^{1} \frac{dP(\lambda z + te^j(z))}{dt} |_{t=0} \; d\lambda \] \[= 
\sum^{n}_{i=1} e^j _i (z)\int^{1}_{0} D_i P(\lambda z)
 d\lambda = \sum^{n}_{i=1} e^j _i (z) T_i(P)(z), \quad 1 \leq j 
\leq n - 1 .\]
Previously, we already noted that $\sum^{n}_{i=1} z_i T_i(P) z = P(z)$. Thus,
the known numbers $T_i(P)(z)$ are solution of the following system of $n$
equations :
\[\left (\begin{array}{ccc}
e^1_1(z) & \mbox{$\ldots$} & \mbox{$e_n ^1(z)$}\\
\vdots & \mbox{$\vdots$} & \mbox{$\vdots$}\\
e^{n-1}_1(z) & \mbox{$\ldots$} & \mbox{$e^{n-1} _n(z)$}\\
z_1 & \mbox{$\ldots$} & \mbox{$z_n$} 
\end{array} \right ) 
\left( \begin{array}{c} T_1(P)(z)\\
\vdots\\
T_{n-1} (P)(z)\\
T_{n}(P)(z)
\end{array} \right) =
\left(\begin{array}{c} 
\int_{\gamma_z}\frac{dP(\lambda z + te^1(z))}{dt} |_{t=0} \; d\lambda \\
\vdots\\
\int_{\gamma_z}\frac{dP(\lambda z + te^{n-1}(z))}{dt} |_{t=0} \; d\lambda \\
P(z) 
\end{array} \right)\]

\noindent The determinant $\Delta (z)$ of the matrix to the left also exists
for $z \in \overline{V_k}$, and it depends
continuously on $z$. It is nowhere zero, and $\overline{V_k}$ is compact, hence
its norm is bounded from below.
The vectors $z$, $e^j (z)\quad (1 \leq j \leq n-1)$ are linearly independent, 
as any complex line passing through $\Omega$ intersects 
$\partial \Omega$ transversally. Hence we can use Cramer's rule to express 
$T_i(P)(z)$ in terms of $\Delta (z)^{-1}$, $e_l^k(z)$, $z$, $P(z)$ and the 
integrals $\int_{\gamma_z}
\frac{dP(\lambda z + te^j(z))}{dt} |_{t=0} \; d\lambda$. Each of those terms 
can be estimated from above with $C ||P||_S$, hence $T_i(P)(z) \leq K_k 
||P||_S$. Since there is only a finite number of $V_j$'s, we have that  
$||T_i(P)||_B \leq K||P||_S$. 
\end{proof}

\section{extending the operators $T_i$}
\noindent Let $f \in H^{\infty} (\Omega)$ such that $f(0)=0$. Let
$p$, $B$, $S$ be as above. Since $\Omega$ is a Runge domain, there is a 
sequence $P_1$, $P_2$, $\ldots$ of polynomials that all vanish at $0$, such 
that $P_n$ converges uniformly to $f$ on $S$. Then $T_i(f)(z) := \lim _{n 
\rightarrow \infty} T_i(P_n)(z)$ exists for all $z \in B$.

\begin{Thm}
The function $T_i(f)$ is properly defined on $\Omega$. Furthermore, it is in 
$H^{\infty}(\Omega)$, and $f(z) = \sum _{i=1} ^{n} z_i T_i(f)(z)$. 
\end{Thm}

\begin{proof}
It is easy to see that the function $T_i(f)$ is properly defined : let $P_n$,
$R_n \rightarrow f$ uniformly on $S$, $T$ respectively. Then
$||T_i(P_n) - T_i(R_n)||_{B} \leq K||P_n - R_n||_{S \cap T} 
\rightarrow 0$.
$||T_i(f)||_B \leq K ||f||_S$, hence $||T_i(f)||_{\Omega} \leq K 
||f||_{\Omega}$. Thus $T_i(f) \in L^{\infty}(\Omega)$. Since the sequence of
polynomials $T_i(P_n)$
converges uniformly on $S$, their limit $T_i(f)$ is in $H(B)$. We also have 
that \[\sum_{i=1}^{n} z_i T_i(f)(z) = \sum_{i=1}^{n} z_i \lim 
_{n \rightarrow \infty} T_i(P_n)(z) = \lim _{n \rightarrow \infty} P_n(z) = 
f(z).\] for $z \in B$. As we can repeat this argument for every point $p \in 
\Omega$ with corresponding neighborhood $B$, the proof is complete.
\end{proof}

\begin{Lem}
Let $f \in A(\Omega)$. For $z \in V_k$ with corresponding $w$,
$w_k \in \partial \Omega$, let $e(z)$ be a complex unit tangent
vector to $\partial \Omega$ at $w$ that varies continuously with $z$. Then 
\[I(z) := \int_{\gamma_{z}} \frac{df(\lambda z + te(z))}{dt}|_{t = 0} \; 
d\lambda \quad \in C(V_k).\]
\end{Lem}

\begin{proof}
This a fairly standard application of the dominated convergence theorem of
Lebesgue. In detail :
let $z \in V_k$, let $V_k \ni z_n \rightarrow z$. Then $\gamma_{z_n}'(s) 
\rightarrow \gamma_{z}'(s)$ for all $s \in [0,1]$. Define $h_n$ and $h$ in 
the following way :
\begin{eqnarray*}h_n(\zeta,s)&:=& \left\{ \begin{array}{ll}
\frac {f(\gamma_{z_n}(s)z_n + \zeta e(z_n))}{\zeta^{2}} 
\gamma^{'}_{z_n}(s),& \mbox{ if  $s \in [0,1- \frac{1}{2} \sigma];$}
\\ \frac {f(\gamma_{z_n}(s)z_n + (1-s)^{1 + \epsilon /2} 
\zeta e(z_n))}{(1-s)^{1 + \epsilon /2} \zeta^{2}} \gamma^{'}_{z_n}
(s), & \mbox{ if $s \in [1 - \frac{1}{2} \sigma,1]$;}\\
\end{array}\right .\\
h(\zeta,s)&:=& \left\{ \begin{array}{ll}
\frac {f(\gamma_{z}(s)z + \zeta e(z))}{\zeta^{2}} 
\gamma^{'}_{z}(s),& \mbox{\qquad if  $s \in [0,1 - \frac{1}{2}
\sigma]$;}
\\ \frac {f(\gamma_{z}(s)z + (1-s)^{1 + \epsilon /2} 
\zeta e(z))}{(1-s)^{1 + \epsilon /2} \zeta^{2}} \gamma^{'}_{z}
(s), & \mbox{\qquad if  $s \in [1 - \frac{1}{2} \sigma,1]$.}\\
\end{array}\right .\end{eqnarray*}
Then
 \begin{eqnarray*}
I(z_n)& =& \int_{\gamma_{z_n}} \frac{df(\lambda z_n + te(z_{n}))}{dt}|_
{t = 0} \; d\lambda  \\
&=&\int _{0} ^{1 - \frac{1}{2} \sigma} \int _{C(0,A)} h_n(\zeta, s)d\zeta ds +
\int _{1 - \frac{1}{2} \sigma} ^{1} \int _{C(0,1)} h_n(\zeta, s)d\zeta ds .
\end{eqnarray*}
For fixed $\zeta$ and $s$, 
$h_n(\zeta,s)$ converges to $h(\zeta, s)$. Furthermore, for all
$n \in \NN$ and $s \in [1 - \frac{1}{2} \sigma,1]$ one has that $|h_n(\zeta, s)| 
\leq M||f||(1- s)^{-1/(1 + \epsilon/2)}$ . 
For $s \in [0,1 - \frac{1}{2} \sigma]$ there is a similar estimate on $h_n$. The function 
$(1- s)^{-1/(1 + \epsilon/2)}$ is integrable on $[0,1]$. Applying 
Lebesgue's theorem yields that \[I(z_n) = \int \int h_n 
\rightarrow \int \int h = I(z),\] thus $I \in C(V_k)$. 
\end{proof}

\noindent Now let $f \in A(\Omega)$, $f(0)=0$. On $\Omega$ we define 
$T_i(f)(z)$ as above. We now proceed to define $T_i(f)$ on $\partial \Omega$.
To every $z \in V_k$ there
correspond $w$, $w_k \in \partial \Omega$. We choose linearly independent 
complex unit tangent vectors $e^1(z)$, $\ldots$, $e^{n-1}(z)$ to 
$\partial \Omega$ at $w$ such that they depend continuously on $z \in V_k$.
Consider the following system of $n$ equations $SY(f)$ :
\[\left (\begin{array}{ccc}
e^1_1(z) & \mbox{$\ldots$} & \mbox{$e_n ^1(z)$}\\
\vdots & \mbox{$\vdots$} & \mbox{$\vdots$}\\
e^{n-1}_1(z) & \mbox{$\ldots$} & \mbox{$e_{n} ^{n-1}(z)$}\\
z_1 & \mbox{$\ldots$} & \mbox{$z_n$} 
\end{array} \right ) 
\left( \begin{array}{c} x_1\\
\vdots\\
x_{n-1}\\
x_n
\end{array} \right) =
\left(\begin{array}{c} 
\int_{\gamma_z}\frac{df(\lambda z + te^1(z))}{dt} |_{t=0} \; d\lambda \\
\vdots\\
\int_{\gamma_z}\frac{df(\lambda z + te^{n-1}(z))}{dt} |_{t=0} \; d\lambda \\
f(z) 
\end{array} \right) \] 
The vectors $z$, $e^j(z) \quad (1 \leq j \leq n-1)$ are linearly independent, 
since any complex line passing through $\Omega$ intersects 
$\partial \Omega$ transversally. Thus the system has a unique solution. 

\begin{Lem}
Let $f \in A(\Omega)$ such that $f(0)=0$. Let $z \in V_{k} ^{\circ}$. Then 
$x_i = T_i(f)(z)$. 
\end{Lem}

\begin{proof}
Choose a compact set $S \subset \Omega$ that contains $z$, as in the proof of 
lemma 2. Choose a compact set $T \subset \Omega$ such that 
$S \subset T^{\circ}$. Let $P_1$, $P_2$, $\ldots$ be a
sequence of polynomials (all vanishing at 0) that converges uniformly to $f$
on $T$. Then $D_i P_n \rightarrow D_i f$ uniformly on $S$, hence, because of 
the chainrule, \[ \lim _{n \rightarrow \infty}
\int _{\gamma_z} \frac{dP_n(\lambda z +te)}{dt} |_{t=0} \; d\lambda =
\int _{\gamma_z} \frac{df(\lambda z +te)}{dt} |_{t=0} \; d\lambda.\] 
In the proof of lemma 2 we saw that for a polynomial $P$, the solution $(x_1$,
$\ldots$, $x_n)$ of the system $SY(P)$ is indeed $(T_1(P)(z)$ ,$\ldots$, 
$T_n(P)(z))$. Taking the limit on both sides of the system of equations 
$SY(P_n)$ yields that $x_i = T_i(f)(z)$.
\end{proof}

\noindent For $z \in \partial \Omega$ we define $T_i(f)(z) := x_i$.

\begin{Lem}
Let $f$ be a function in $A(\Omega)$ that vanishes at $0$. Then $T_i(f) \in
A(\Omega)$.
\end{Lem}

\begin{proof}
Let $\Delta (z)$ be the determinant of the matrix to the left in $SY(f)$. 
We again use Cramer's rule to express $T_i(f)(z)$ in terms of 
$\Delta (z)^{-1}$, $e_l^k(z)$, $z$, $f(z)$ and the integrals 
$\int_{\gamma_z} \frac{df(\lambda z + te^j(z))}{dt} |_{t=0} \; d\lambda$. 
These are all continuous functions of $z$ on $V_k$. Therefore $T_i(f)$ is in 
$C(V_k)$, and repeating this argument for all $k$ yields that
$T_i(f) \in C(V)$. Hence $T_i(f) \in A(\Omega)$.
\end{proof}

\begin{Thm}
Let $\Omega \subset \CC^n$ be a linearly convex domain with $C^{1 + \epsilon}$-
boundary. Then the ideal in $A(\Omega)$ (or $H^{\infty}(\Omega)$) 
consisting of all functions in $A(\Omega)$ (or $H^{\infty}(\Omega)$) that 
vanish at $p \in \Omega$ is generated by the coordinate functions $z_1 - p_1,
\ldots, z_n - p_n$.
\end{Thm}

\begin{proof}
For an $f \in A(\Omega)$ (or $H^{\infty}(\Omega)$), such that $f(0)=0$, we have
that \[f(z) = \sum _{i=1} ^{n} (z_i - p_i) T_i(f)(z),\] and $T_i(f) \in 
A(\Omega)$ (or $H^{\infty}(\Omega)$).
\end{proof}

\vskip5mm

\section{Final remarks.}
\label{Final remarks.}  

\noindent
In \cite{gran} Grang\'{e} gave the following example
of a convex domain in $\CC ^2$ for which $T_2(f)$ is unbounded for a certain
$f \in H^{\infty}(\Omega)$ :
let $h(x):= \frac{-x}{\log x}$ for $x>0$, $h(0):=0$. 
Let \[\Omega:=\{(z_1,z_2) \in 
\CC^2 : |z_2| < 1, |z_1|^{2} + h(|z_2|) - 1 < 0 \}.\]
This shows that the functions $T_i(f)$ may fail to solve the Gleason problem
on \\$\CC$-convex domains with $C^1$-boundary.
However, it is possible to solve the Gleason problem for $H^{\infty}(\Omega)$
and $A(\Omega)$ by using different techniques, as the present authors showed
in \cite{LeWi}.

\vskip5mm \noindent
A glance at the previous proofs may suggest that our results can be obtained
under the weaker assumptions of the next lemma. This lemma however shows
that these assumptions are not really weaker at all.

\begin{Lem}
Let $\Omega$ be a bounded domain with $C^1$-boundary such that every complex
line passing through $\Omega$ intersects $\partial \Omega$ transversally.
Suppose that $\Omega$ intersected with any complex line is connected. Then
$\Omega$ is $\CC$-convex.
\end{Lem}

\begin{proof}
From the conditions it follows immediately that $\Omega$ is connected. Suppose
$\Omega$ is not $\CC$-convex. Then it is not weakly linearly convex either,
meaning there is a point $z \in \partial \Omega$ such that every complex
hyperplane $H$ through $z$ intersects $\Omega$. We take for $H$ the complex 
tangent space to $\partial \Omega$ at $z$. It contains a complex line that is 
tangential to $\partial \Omega$ at $z$ and intersects $\Omega$. This 
contradicts our assumption that such a line intersects $\partial \Omega$
transversally.
\end{proof}

\vskip10mm \noindent
Department of Mathematics\\
University of Amsterdam\\
Plantage Muidergracht 24\\
1018 TV Amsterdam\\
The Netherlands


\vskip10mm

\begin{thebibliography}{15}

\bibitem{APS1} Andersson, M., Passare, M. and R. Sigurdsson, {\em Complex
convexity and analytic functionals I}, Report RH-06-95, Science Institute, 
University of Iceland, 1995.
\bibitem{APS2} Andersson, M., Passare, M. and R. Sigurdsson, {\em Complex
convexity and analytic functionals II}, Report RH-20-2000, Science Institute, 
University of Iceland, 2000.
\bibitem{bafa1} Backlund, U. and A. F\"{a}llstr\"{o}m, {\em On Gleason's
problem for $H^{\infty}$}, Research Report no.{\bf 6} (1992), Department of 
Mathematics, University of Ume\aa, Sweden.
\bibitem{bafa2} Backlund, U. and A. F\"{a}llstr\"{o}m, {\em The Gleason problem
for $A(\Omega)$}, New Zealand J. Math. {\bf 24} (1995), no.1, 17--22.
\bibitem{bafa3} Backlund, U. and A. F\"{a}llstr\"{o}m, {\em Counterexamples to
the Gleason problem}, Ann. Scuola Norm. Sup. Pisa Cl. Sci. (4) {\bf 26} (1998),
no.3, 595--603.
\bibitem{bafa4} Backlund, U. and A. F\"{a}llstr\"{o}m, {\em The Gleason 
property for Reinhardt Domains}, Math. Ann. {\bf 308} (1997), 85--91.
\bibitem{foev} Forn\ae ss, J. E. and N. \O vrelid, {\em Finitely generated 
ideals in $A(\Omega)$}, Ann. Inst. Fourier (Grenoble) {\bf 33} no. 2 
(1983), 77--85.
\bibitem{glea} Gleason, A. M., {\em Finitely generated ideals in Banach
 algebras}, J. Math. Mech. {\bf 13} (1964), 125--132.
\bibitem{gran} Grang\'{e}, M., {\em Diviseurs de Leibenson et probl\`{e}me de 
Gleason pour $H^{\infty}(\Omega)$ dans le cas convexe}, Bull. Soc. Math. 
France {\bf 114} (1986), 224--245.
\bibitem{Hor} H\"{o}rmander, L., {\em Notions of convexity}, Birkh\"{a}user
Boston, Inc., Boston, MA, 1994.
\bibitem{khen} Khenkin, G. M., {\em Approksimatsiya funktsij v psevdovypuklych
oblastyach i teorema Z. L. Leibenzona}, Bull. Acad. Polon. Sci. S\'{e}r. Sci. 
Math. Astronom. Phys. {\bf 19} (1971), 37--42. (Russian)
\bibitem{LeWi} Lemmers, O. and J. Wiegerinck, {\em Reinhardt domains and the 
Gleason problem}, to appear in Ann. Scuola Norm. Sup. Pisa Cl. Sci.
\bibitem{noel} Noell, A., {\em The Gleason problem for domains of 
finite type}, Complex Variables Theory Appl. {\bf 4} (1985), 233--241.
\bibitem{oev1} \O vrelid, N. {\em Integral representation formulas and
 $L^p$ estimates for the $\bar{\partial}$-equation}, Math. Scand. 
{\bf 29} (1971), 137--160.
\bibitem{oev2} \O vrelid, N. {\em Generators of the maximal ideals of 
$A(\overline{D})$}, Pacific J. Math. {\bf 39} (1971), 219--223.

\end{thebibliography}
\end{document}